\documentclass[twocolumn]{autart}    
\usepackage{xcolor}

\usepackage{graphicx}          

\usepackage{amsmath}
\usepackage{amssymb}
\usepackage{mathtools}
\usepackage{enumitem}
\usepackage{comment}
\usepackage{todonotes}
\newcommand{\Real}{\mathbb{R}}
\newcommand{\Neal}{\mathbb{N}}
\newcommand{\Peal}{\mathbb{P}}
\newcommand{\Eeal}{\mathbb{E}}

\newcommand{\calX}{\mathcal{X}}
\newcommand{\calV}{\mathcal{V}}
\newcommand{\calI}{\mathcal{I}}
\newcommand{\calU}{\mathcal{U}}

\newcommand\calM{\mathcal{M}}

\DeclareMathOperator{\Tr}{tr}
\DeclareMathOperator{\dom}{dom}

\begin{document}

\begin{frontmatter}

\title{Characterization of Safety in Stochastic Difference Inclusions using Barrier Functions\thanksref{footnoteinfo}} 

\thanks[footnoteinfo]{Will be added in subsequent versions. }

\author[Boulder]{Masoumeh Ghanbarpour}\ead{masoumeh.ghanbarpour@colorado.edu},    
\author[Boulder]{Sriram Sankaranarayanan}\ead{srirams.@colorado.edu}

\address[Boulder]{Department of Computer Science, University of Colorado Boulder, USA.}  


\begin{abstract}                          
 We study stochastic systems characterized by difference inclusions. Such stochastic differential inclusions are  defined by set-valued maps involving the current state and stochastic input. For such systems, we investigate the problem of proving bounds on the worst-case probability of violating safety properties. Our approach uses the well-known concept of barrier functions from the study of stochastic control systems. However, barrier functions are hard to prove in the presence of stochastic inputs and adversarial choices due to the set-valued nature of the dynamics. 
 
 In this paper, we show that under some  assumptions on the set-valued map including upper semi-continuity and convexity combined with a concave barrier function vastly simplifies the proof of barrier conditions, allowing us to effectively substitute each random input in terms of its expectation. We prove key results based on the theory of set-valued maps and provide some interesting numerical examples. The ideas proposed here will contribute to the growing interest in problems of robust control and verification of stochastic systems in the presence of uncertain distributions and unmodeled dynamics. 
\end{abstract}

\end{frontmatter}

\section{Introduction}
Ensuring the safety of dynamical systems is particularly challenging when their dynamics are influenced by uncertainty. This uncertainty may arise from multiple sources with varying mathematical properties and implications on the system
behavior. Common types include stochastic uncertainties modeled by random variables with known probability distributions
and set-valued uncertainties modeled by a \emph{non-deterministic} (adversarial) choice from a set of possibilities. 
In this paper, we consider the analysis of discrete-time dynamical systems that combine uncertainties from
two distinct sources: stochastic and set valued uncertainties. We express the dynamics as $x^+ \in F(x, v)$, wherein $x$ is the current state of the system, $x^+$ is the next state, $v$ is a random input drawn from a fixed probability distribution and $F$ is a set valued map. Since $F$ is set valued, we assume that an ``adversary'' selects a specific next state in the set $F(x, v)$.
The combination of two sources allows us to model interesting problems that include (a) systems whose dynamics are influenced by stochastic and set-valued disturbances; (b) systems with unmodeled dynamics that may be incorporated into the set-valued maps; and (c) systems with stochastic inputs drawn from a set (family) of distributions rather than a single fixed distributions wherein the distributions may even be time-varying.

Our goal is to analyze the probability of violating safety properties under these combined uncertainties.  Specifically, we investigate the problem of proving bounds on the probability that a system governed by such a set-valued map will enter a given set of unsafe states $X_u$ over some finite or infinite time horizon of interest. In this case, we will assume that the adversarial choices will be made at each step to ensure that the system will enter $X_u$ within the given time horizon. To do so, we extend the classic notion of a barrier function first introduced for deterministic systems~\cite{prajna2004safety}, and then extended to stochastic systems~\cite{prajna2004stochastic,steinhardt2012finite}. 

In this paper, we first make assumptions on the set-valued map that will allow us to formally define notions
such as the ``worst-case'' probability of reaching a set $X_u$. We demonstrate that if $F$ is chosen to be a bounded and upper semicontinuous set valued map with key measurability properties (Cf. Assumption~\ref{ass:standing}), the worst-case probabilities are well-defined. Likewise, we show that for nonnegative and upper semicontinuous barrier functions, the worst-case expected value at a future time is well-defined and can be used to compute  upper bounds on the worst-case probabilities of reaching certain sets in the state-space. As a result, barrier functions provide a general framework for verifying safety. However, the expectations of the barrier functions can be hard to compute, in practice since they involve a supremum over possible choices of the next states from a set of possibilities. To enable easier computation of the bounds, we study the case where the barrier function is assumed to be concave while the set-valued map is assumed convex. We demonstrate that this additional structure leads to stronger and more easily verifiable sufficient conditions for both safety and viability.
Following this, we introduce the concept of a supermartingale barrier functions. Supermartingale properties naturally capture the tendency of the barrier function to decrease in expectation, which is well-suited to ensuring safety in stochastic systems. We also show that the concavity assumption with supermartingale leads to a more applicable supermartingale condition that is easy to verify and yields bounds that are valid for an infinite time horizon. The paper demonstrates the key ideas through some numerical results obtained on interesting classes of set-valued stochastic systems.

\subsection{Related Work}
The analysis of safety and reachability properties of probabilistic systems has been well-studied. Abate et al~\cite{abate2008probabilistic} investigate probabilistic reachability for discrete-time stochastic hybrid systems, employing dynamic programming to compute reachable sets and establish probabilistic safety guarantees. In this paper, we borrow the terminology of \( p \)-safety, introduced by Wisniewski et al~\cite{wisniewski2021safety}, which defines a system as $p$-safe if  the probability of the system avoiding an unsafe set is at least \( 1 - p \). This method leverages sum-of-squares techniques to systematically quantify and compute safety probabilities. Recent works, such as Chern et al~\cite{chern2021safe} and Lavaei et al~\cite{lavaei2022safety}, focus on synthesizing control laws to ensure the safety of stochastic dynamical systems.

The area of distributionally robust optimization studies optimization problems with chance constraints wherein the distributions governing the random variables belong to a family of possible distributions~\cite{Rahimian+Mehrotra/2022/Frameworks}. This area has incorporated interesting approaches to model families of distributions including distribution parameters/moments drawn from a set~\cite{Bertsimas+Popescu/2005/Optimal,Delage2010Distributionally} and ideas from optimal transport, including families defined by bounded Wasserstein distances~\cite{Pflug+Wozabal/2007/Ambiguity}. More recently, there has been a lot of interest on the robust control of stochastic systems, wherein the distributions of the stochastic inputs are drawn from a family rather than a fixed probability distribution~\cite{Schon+Others/2024/Data,Taha+Others/2023/Distributionally}. Our approach using set-valued maps can also capture some forms of distributional robustness in the model of the system. For instance, if the distributions range over families defined by a set of possible values of the distribution parameters/moments, it is possible in some cases to express the resulting system as a set-valued map using randomized inputs drawn from a fixed distribution.

Several studies explore the interplay between stochastic disturbances and adversarial inputs. Summers et al \cite{summers2010verification} propose a stochastic games framework for verification and control of discrete-time stochastic hybrid systems, addressing scenarios where an adversary influences system dynamics with objectives counter to those of the controller. In the presence of distributional uncertainties, Yang et al \cite{yang2018dynamic} introduce a dynamic game approach to distributionally robust safety specifications, assuming that the disturbance distribution is only partially known and belongs to an ambiguity set of probability distributions. Furthermore, Chapman et al~\cite{chapman2021risk} employ Conditional Value-at-Risk (CVaR) for risk-sensitive safety analysis, capturing rare but harmful outcomes in stochastic systems and quantifying expected losses in worst-case scenarios. The concept of probabilistic invariance has also been explored for enforcing safety in stochastic systems. Methods such as linear programming and mixed-integer linear programming have been utilized to compute invariant sets that ensure safety over both finite and infinite horizons \cite{schmid2023probabilistic}.

The problem of reasoning about combinations of ``demonic'' nondeterminism and probabilistic inputs has been studied by computer scientists
seeking to prove properties of programs with random number generating statements and set-valued uncertainties~\cite{McIver2005,Batz+Others/2024/Programmatic}. The reasoning is performed by a systematically propagating logical assertions about the state of the program upon termination backwards through the weakest pre-expectation operator until an assertion over the known initial states of the program is verified. While this approach can be quite general and subsume some of the ideas here, significant modeling and computational difficulties  arise when the program involves random variables with an infinite set of support, or non-determinism that can make infinitely many choices. In contrast, the approach taken by systems and control theorists (including the approach of this paper) focuses on systems that work over real-valued system variables and random variables that go beyond finitely many choices. At the same time, the system structure is taken to be a map whose form remains unchanged rather than a computer program that has nontrival control flow and may contain simple/nested loops, function calls and recursion. 

The stability of stochastic difference inclusions has also been well studied~\cite{possieri2017lyapunov}. Maghenem et al derive sufficient conditions for robust safety in differential inclusions using barrier functions, allowing for unbounded safety regions and nonsmooth barrier functions~ \cite{maghenem2025sufficient}. More recently Maghenem and Ghanbarpour have established the equivalence between robust safety and the existence of a barrier function certificate for differential inclusions, constructing a barrier function based on the time-to-impact function with respect to a specifically designed reachable set~\cite{maghenem2025converse}. The key difference here is that, to our knowledge, there has been no work done on the combination of set-valued maps with probabilistic inputs, which is the main topic of this paper.

\textbf{Notation.} 
For sets \( K,O \subseteq \mathbb{R}^n \), \( K \backslash O \) denotes the subset of elements of \( K \) that are not in \( O \). $\mathbb{B}$ denotes the open unit ball centered at the origin. $U(x)$ indicates a neighborhood set around the point $x$.
For a function \( f : \mathbb{R}^n \rightarrow \mathbb{R}^m \), \( \operatorname{dom} f \) denotes the domain of definition of \( f \), and if \( f \) is differentiable, \( \nabla f \) denotes the gradient of \( f \). \( \overline{\mathbb{N}} \) represents \( \mathbb{N} \cup \{\infty\} \).
A set-valued map \( F : \mathbb{R}^m \rightrightarrows \mathbb{R}^n \) associates each element \( x \in \mathbb{R}^m \) with a subset \( F(x) \subseteq \mathbb{R}^n \).
$S^n_{++}$ denotes $n \times n$ matrices which are symmetric positive definite. For matrix $X$, $\Tr(X)$ indicates trace of $X$. $\vee$ denotes or, and $[1,N]:=\{1,2,\cdots,N\}$.
Given \( C \subseteq \mathbb{R}^n \), \( \mathcal{I}_C:\mathbb{R}^n \to \{0,1\} \) is the indicator function for the set \( C \), defined as
$\mathcal{I}_C(x) :=
\begin{cases}
1, & \text{if } x \in C, \\
0, & \text{otherwise}.
\end{cases}$
\section{Preliminaries}
In this section, we present fundamental concepts of continuity, followed by the definition of a stochastic difference inclusion system. We then state the notions of reachability and safety in a probabilistic context.

\subsection{Continuity Notions}
We recall some of the continuity notions for set-valued and single-valued maps.
Let $F: \calX \rightrightarrows \mathbb{R}^m$ be a set-valued map, where $\calX \subseteq \mathbb{R}^n$.

\begin{defn} [Upper Semicontinuous Maps]
A set-valued map $F$ is  \textit{upper semicontinuous} at $x \in \calX$ if for each 
$\delta > 0$, there exists a neighborhood $U(x)$ of $x$,  such that for each $y \in U(x) \cap \calX$, $F(y) \subseteq F(x) + \delta \mathbb{B}$; see \cite[Definition 1.4.1]{aubin2009set}.
\end{defn}
\begin{defn}[Locally Bounded] \label{def:locallybdd}
 A set-valued map  $F$ is \textit{locally bounded} at $x \in \calX$,  if there exists a neighborhood  $U(x)$ of $x$, and a constant $\beta(x) > 0$ such that for all  $y \in U(x) \cap \calX$ and 
  $\zeta \in F(y)$,  $\| \zeta \|  \leq \beta(x)$.  
\end{defn}
The map $F$ is said to be upper semicontinuous if it is upper semicontinuous for all $x\in \calX$.
We will recall the definition of upper semicontinuity for functions (singleton set-valued maps).
\begin{defn}[Upper Semicontinuous Functions] \label{def:cont_sing}
A scalar function $B: K \rightarrow \mathbb{R}$, where $K \subseteq \mathbb{R}^n$ is said to be \textit{upper semicontinuous} at $x \in K$ if, for every sequence $\left\{ x_i \right\}_{i=0}^{\infty} \subseteq K$ such that $\lim_{i \rightarrow \infty} x_i = x$, we have $\limsup_{i \rightarrow \infty} B(x_i) \leq B(x)$.

The function \( B \) is said to be \emph{upper semicontinuous} if it satisfies this property for all \( x \in K \). Moreover, if \( B \) is upper semicontinuous, then \( -B \) is \emph{lower semicontinuous}.

\end{defn} 
\subsection{Difference inclusion, Reachability, and Viability} \label{subsec:dinc_reach}

Consider a difference inclusion written as: 
\begin{equation} \label{eq:sys_sp}
    x^+ \in F(x), \quad x \in \calX\subseteq \Real^n,
\end{equation}
where $x \in \calX$ is the state, $ F: \calX \rightrightarrows \mathbb{R}^n $ is a set-valued map, and $x^{+}$ denotes state in the next time step. 

We view set-valued maps in terms of a game played between a system and an adversary.
For each $x$, the set \( F(x)  \) represents all the possible choices available to the adversary for the next time step. I.e, the adversary selects an element \( \hat{x} \in F(x) \).

A set \( S \subseteq \mathbb{R}^n \) is said to be \textit{reachable} from \( x \) in a single time step if \( F(x) \cap S \neq \emptyset \). This means that the adversary can choose an element \( \hat{x} \in F(x) \) such that \( \hat{x} \in S \), ensuring that the next state lies within \( S \). Conversely, if \( F(x) \cap S = \emptyset \), then no choice from \( F(x) \) allows the system to reach \( S \) in one step.

Similarly, the solution from $x$  is \textit{viable} for a single time step in the set \( S \) if \( F(x) \subseteq S \) for all \( x \in S \). In this case, regardless of the adversary’s choice, the system remains within \( S \) in the next time step. However, if there exists some \( x \in S \) for which \( F(x) \not\subseteq S \), then the adversary can choose an element that leads the system out of \( S \), violating viability.

\subsection{Stochastic Difference Inclusion}  
A stochastic difference inclusion system is a set-valued map with  a stochastic input drawn from a distribution and is given as
\begin{equation} \label{eq:sys}
    \Sigma: \quad x^{+} \in G(x,v),  \quad x \in \calX \subseteq \mathbb{R}^n, \, v \in \calV\subseteq \Real^m, 
\end{equation}
where $G:\calX \times \calV \rightrightarrows \mathbb{R}^n$ is the set-valued map, $x$ is the state, and $v$ is a random input with a measure $\mu$  over a the support $\calV$.
Here, we study the safety in probability for the system $\Sigma$.

The key assumptions for the system and the solution are given in the following.

\begin{assum} \label{ass:standing}
The set-valued mapping $G$ in \eqref{eq:sys} satisfies the following conditions:
\begin{enumerate}
    \item \label{item:C1}The map $x \mapsto G(x,v)$ is upper semicontinuous for each $v \in \calV.$
    \item \label{item:C2} The map $G$ is locally bounded (see Definition \ref{def:locallybdd}).
    \item For each $(x,v)$, $G(x,v)$ is nonempty, closed, and convex.
    \item \label{item:C3} The mapping $v \mapsto \operatorname{graph}(G(.,v))=\{(x,y)\in \calX \times \calX\,:\, y \in G(x,v)\}$ is measurable. 
\end{enumerate}
\end{assum}

\begin{assum} \label{ass:iid}
Let $(\Omega, \mathcal{F}, \mathbb{P})$ be a probability space \cite{folland1999real}.  
Consider a sequence of random vectors $(v_i)_{i \in \mathbb{N}}$, where each  
$v_i : \Omega \to \mathcal{V} \subseteq \mathbb{R}^m$.  
We assume that $(v_i)_{i \in \mathbb{N}}$ are independent and identically distributed (i.i.d.) with common distribution $\mu$ on $(\mathcal{V}, \mathcal{B}(\mathcal{V}))$, where $\mathcal{B}(\mathcal{V})$ denotes the Borel $\sigma$-algebra on $\mathcal{V}$.
\end{assum}

\begin{defn}[Trajectory]\label{def:trajectory}
Let $v_0, \ldots, v_{k-1}$ be a finite sequence of samples drawn i.i.d. from the distribution $\mu$.
 A \emph{trajectory}
of the system starting from state $x_0 \in \calX$, and corresponding to the sequence $v_0, \ldots, v_{k-1}$ is a sequence of states $x_0, \ldots, x_{k}$ obtained as follows:
\begin{equation}\label{eq:set-valued-map-iteration}
\left\{ \begin{array}{rll}
x_0 & = x_0 \\ 
X_{i+1} & = G(x_i, v_i) &  i = 0, \ldots, k-1, \\ 
x_{i+1} & \in X_{i+1} & i = 0, \ldots, k-1\,.\\ 
\end{array} \right.
\end{equation}
We will use the notation $X_0, \ldots, X_k$ to denote the sequence of sets that correspond to the solution $x_0, \ldots, x_k$, wherein $X_0 = \{ x_0 \}$.
\end{defn}

Note that since $G(x_i, v_i)$ is nonempty (Cf. assumption \ref{ass:standing}), we can obtain trajectories of any length $k \in \mathbb{N}$ starting from any $x_0 \in \calX$. Also, note that multiple trajectories of length $k+1$ may exist for a fixed sequence of samples $v_0, \ldots, v_{k}$.

In the following Lemma we show that $G(x_i, v_i)$ is measurable for each $i\geq 0$.
\begin{lem}
    [Measurability] \label{lm:meas}
    Consider the system $\Sigma$ described by \eqref{eq:sys}, and suppose Assumption \ref{ass:standing} holds. Then, for each $k\geq 0$, $G(\cdots (G(G(x,v_0),v_1),v_2),\cdots),v_k)$ is closed and measurable for all $x$.
\end{lem}

\begin{pf}
    See Appendix \ref{pf:lm:meas} for the proof.
\end{pf}

We will now present a few illustrative examples of models based on the formalism presented thus far.

\begin{exmp}[Set-Valued Random Walks]
Consider a state $x \in \Real$, modeling the aggregate water level in a large tank that can be filled up by rain water. The total rainfall each month is modeled 
as a normal random variable with a time-varying mean $\mu(t)$ and time-varying 
standard deviation $\sigma(t)$. Furthermore, the system is subject to an uncertain monthly demand $w(t)$ which belongs to the 
set $w(t) \in [a, b]$. 
\begin{equation}
  x(t+1) \in x + \mu(t) + \sigma(t) v - [a, b] \,, v \sim N(0, 1)\,.
\end{equation}
In this example, the set-valued map formalism combines two sources of uncertainty: one that is stochastic in nature and the other that is nondeterministic. Suppose we were to bound the probability that the tank starting from a level of $x(0) = 10$, empties out within $6$ months, or overflows. Rather than assume worst-case rainfall, we would like to treat the rainfall as a stochastic input while treating demand as an adversarial worst-case. For instance, when analyzing the probability of emptying the tank, it is natural to assume that the demand is always $b$ during the time period under consideration. Alternatively, when analyzing the probability of the tank overflowing, the demand is set to $a$, which happens to be the worst case for that scenario. In the presence of more complicated dynamics involving multiple state variables and non-linear maps, the analysis can be more involved. 
\end{exmp}

\begin{exmp}[Distributional Robustness]
Consider a stochastic system $x^+ = g(x, v)$ wherein $v$ is a random variable whose distribution is drawn from a family of distributions. For instance, $v$ may be drawn from a set of normally distributed random variables with mean $\mu \in [\mu_l, \mu_u]$ and  standard deviation $\sigma \in [\sigma_l, \sigma_u]$. In this instance, we may write $v = \mu + \sigma \hat{v}$, wherein 
$\hat{v}$ is a standard normal random variable with mean $0$ and standard deviation $1$. In fact, many families of random variables have such convenient parameterizations wherein each instance of the family can be modeled as a function $v = f(\theta, \hat{v})$, wherein $\hat{v}$ is a standard random variable and $\theta$ represents parameters that define the family~\cite{johnson1994continuous}. 

In our framework, if $\theta \in \Theta$ for a set of possible parameters $\Theta$, we can incorporate our system into a set valued map of the form  $x^+ \in G(x, \hat{v})$, wherein,  
\[G(x, \hat{v}) = \{ \hat{x}\ |\ \exists \theta \in \Theta,\ \hat{x} = g(x, f(\theta, \hat{v})) \} \,. \]
Note that the dynamics defined by the  map $G$ is more general in that it can allow the actual distribution parameter $\theta$ to be time-varying rather than fixed. Modeling $\theta$ as a time invariant choice can be achieved if it were incorporated as a state variable in the system description. 
\end{exmp}

\subsection{Probabilistic Reachability and Safety}
First, we define reaching a set within $k$ time steps.
\begin{defn}[Reachability]\label{def:reach-seq}
A nonempty set $S \subseteq \mathcal{X}$ is said to be \emph{reachable within $k \geq 1$ steps} from an initial state $x_0 \in \mathcal{X}$ if there exists a solution of $\Sigma$ starting at $x_0$ that enters the set $S$ in at most $k$ steps.
\end{defn}

Consider the system $\Sigma$ in \eqref{eq:sys}. At each step in~\eqref{eq:set-valued-map-iteration}, we iterate between application of the map $G$ yielding a set $X_{i+1}$ and the choice of a state $x_{i+1} \in X_i$.  This choice is made by an ``adversary'' whose goal, for the purposes of this paper, is to force the system to enter a designated set of states $S$. 

We define a  \emph{utility function} for an adversary as 
$u(x) = \mathcal{I}_S(x)$,
where $\mathcal{I}_S$ is the indicator function of $S$. In other words, at any state $x$, the adversary derives an utility of $1$ if $x \in S$ and $0$ otherwise. 

\paragraph*{Single Step Reach Probability:} With this definition of utility, we can define  the probability of reaching the set $S$ in one step starting from state $x$ as the expectation over the maximum possible utility for an adversary in the next step.

\begin{align} \label{eq:cap1}
     \Peal( G(x,v_0) \cap S\not= \emptyset|x) = \Eeal\left[ \sup_{y \in G(x,v_0)} \calI_S(y)\right]\,.
\end{align}
The lemma below proves the existence of this probability for one-step reachability.
\begin{lem} \label{lm:meas_well}
   Given system $\Sigma$ described by \eqref{eq:sys} satisfying Assmpt. \ref{ass:standing}--\ref{ass:iid}, let $h:\calX \to \Real_{\geq 0}$ be a upper semicontinuous function, then  $ \lambda(x,v):=\displaystyle\sup_{y \in G(x,v)} h(y)$ is measurable. 
   
\end{lem}
\begin{pf}
    See Appendix \ref{pf:lm:meas_well} for the proof.
\end{pf}

As a result, the probability that a given set $S$ is reachable in \emph{one step} is well-defined, since the expectation $\Eeal\left[ \displaystyle\sup_{y \in G(x,v_0)} \calI_S(y) \right]$, wherein the function $\calI_S(y)$ is upper-semicontinuous.  Furthermore, the adversary encounters a set $X_1 = G(x, v_0)$ and either chooses a state $x_1 \in X_1 \cap S$, or fails to do so if $X_1 \cap S = \emptyset$.

\paragraph*{$N$-step reachability:} We define the sequence of functions $\{\lambda_k^S\}$, which $\lambda_k^S:\calX \times \calV^k\to \Real_{\geq 0}$ indicates the supremum of the indicator function at the step $k$. 
\begin{equation}  \label{eq:lambda_S}
\left. \begin{array}{rl}
\lambda_0^S(x) & = \calI_S(x) \\ 
\lambda_1^S(x, v_0) & = \displaystyle\sup_{x_1 \in G(x, v_0)} \lambda^S_0(x_1) \\ 
\lambda_2^S(x, v_0, v_1) & = \displaystyle\sup_{x_1 \in G(x, v_0)} \lambda_1^S(x_1, v_1) \\ 
& \vdots \\ 
\lambda_{k+1}^S(x, v_0, \ldots, v_{k}) & = \displaystyle\sup_{x_1 \in G(x, v_0)} \lambda_{k}^S(x_1, v_1, \ldots, v_{k}) \\ 
\end{array}\right\}
\end{equation}
Notice that in the function $\lambda_{k+1}$ supremum is taking from all the regions that the solution to $\Sigma$ can reach at time step $k+1$ starting from $x$, as
\begin{align} \label{lambda_supsup}
    \lambda_{k+1}^S&(x, v_0, \ldots, v_{k})= 
    \\&\sup_{y_1 \in G(x,v_0)} 
    \sup_{y_2 \in G(y_1,v_1)} \cdots \sup_{y_{k+1} \in G(y_{k},v_k)} \calI_S(y_{k+1}). \nonumber
\end{align}

\begin{lem} \label{lm:lambda_meas_k}
Consider the system $\Sigma$ described by \eqref{eq:sys} such that Assumptions \ref{ass:standing}--\ref{ass:iid} hold.
Let $S\subset \calX$ be a closed set. For each $k \in \mathbb{N}$ and for all $x\in \calX$, the function $\lambda_{k+1}^S(x, v_0, \ldots, v_k)$, defined in \eqref{eq:lambda_S}, is measurable.
\end{lem}

\begin{pf}
    See Appendix \ref{pf:lm:lambda_meas_k} for the proof.
\end{pf}

\begin{lem}
Consider the system $\Sigma$ in \eqref{eq:sys} such that Assumptions \ref{ass:standing}--\ref{ass:iid} hold.
Let $S \subset \calX$ be a closed set, then the probability of reaching $S$ within $N$ time steps from $x \in \calX$ can be bounded as 
\begin{align} \label{lambda_N_I}
    \Peal \big(\exists k\in [1,N], \lambda_k^S &(x,v_0,\cdots,v_{k-1})\geq 1\mid x \big) \\&
    \leq \Eeal\big[ \sup_{k\in [1,N]} \lambda_k^S(x,v_0,\cdots,v_{k-1})\mid x  \big]. \nonumber
\end{align}
\end{lem}
\begin{pf}
  Using the equation \eqref{lambda_supsup}, $\lambda_{k}$ is either $1$, or $0$. If there is possibility that a solution starting from $x$ reaches to the set $S$ at step $k$, $\lambda_k$ is one otherwise it is zero.
Therefore, the probability of reaching the set $S$ within $N$ steps is given as
\begin{align*} 
    \Peal \big(\exists k\in [1,N], &\lambda_k^S(x,v_0,\cdots,v_{k-1})\geq 1\mid x \big) \\&=
    \Peal \big(\sup_{k\in [1,N]} \lambda_k^S(x,v_0,\cdots,v_{k-1})\geq 1\mid x \big) \nonumber\\&
    \leq \Eeal\big[ \sup_{k\in [1,N]} \lambda_k^S(x,v_0,\cdots,v_{k-1})\mid x  \big]. \nonumber
\end{align*}  
\end{pf}

Probabilistic safety is  defined by bounding the probability of reaching the unsafe set starting from an initial set.
\begin{defn} \label{def:safety_prob}
(Probabilistic Safety)  
The system $\Sigma$ is safe \textbf{within} $N \in \overline{\mathbb{N}}$ time steps with probability $1-\rho$, wherein $\rho \in [0,1)$, with respect to the sets $(X_o, X_u)$, where $X_o \subseteq \calX$ is the initial set and $X_u \subseteq \calX$ is the unsafe set with $X_o \cap X_u = \emptyset$, if for any $x \in X_o$,  
\begin{equation} \label{eq:safety_prob_general}
    \Peal( \exists k \in [1, N], \lambda^{X_u}_{k}(x, v_0, \ldots, v_{k-1}) \geq 1 )< \rho \,.
\end{equation}
\end{defn}


\section{Barrier Functions}

 In the previous section, we use the indicator function $\calI_S$ to find a probability of reaching the set $S$. 
 In certifying the safety we wish to bound the probability of reaching the unsafe set. 
 Using indicator function, for bounding the probability of reaching the unsafe set by $\rho$ from $x$ within $N$ times step, we have the following condition
 \begin{equation} \label{eq:reach_S_I}
     \Eeal\big[ \sup_{k\in [1,N]} \lambda_k^S(x,v_0,\cdots,v_{k-1})\mid x  \big] \leq \rho. 
 \end{equation}
 In this section, we study the use of other nonnegative functions instead of the indicator function. 
Using such functions may provide better conditions than the one in~\eqref{eq:reach_S_I}, 
since we consider the unsafe set as a $\Delta$-superlevel set of these functions, 
where $\Delta > 0$. Consequently, the right-hand side of the condition changes to $\rho \Delta$.

 We show how to modify the indicator function using a Barrier function to obtain bounds on the probability of reaching an unsafe set. 
 Consider a nonnegative upper semicontinuous function $B:\calX \to \Real_{\geq 0}$. 
The $\Delta$-sublevel set of $B$ for $\Delta\geq 0$ is defined as
\begin{equation} \label{eq:k}
    K_\Delta:= \{x \in \calX\,:\, B(x) \leq \Delta\}\,.
\end{equation}
If $B$ is upper semicontinuous, then every sublevel set of $B$ is closed. Therefore, $K_\Delta$ is a closed set. Let $K_{\Delta}^C$ denote the
set  $\mathcal{X} \setminus K_\Delta$.
The condition of one step reachability from $x$ to $K_{\Delta}^C$ in $\calX$ is \( G(x,v) \cap K_{\Delta}^C \neq \emptyset \) for a given $v \in \calV$, and can 
be equivalently expressed as $\displaystyle\sup_{x^+ \in G(x,v)} B(x^+)> \Delta$.


Starting from $x\in \calX$, the probability of reaching $K_{\Delta}^C$ in the next time step using the barrier function is given as
\small{\begin{align} \label{eq:test}
    \mathbb{P}( G(x,v) \cap K_{\Delta}^C \neq \emptyset \mid x)&=\mathbb{P} \left( \displaystyle \sup_{x^{+}\in G(x,v)} B(x^{+})  > \Delta \mid x\, \right)
\end{align}}

Continuing with the equation \eqref{eq:test}, and using Markov inequality, we have
\small{\begin{equation}
    \mathbb{P} \left( \displaystyle \sup_{x^+ \in G(x,v)} B(x^+) > \Delta \,\mid x\, \right) \leq \frac{1}{\Delta} \Eeal[\sup_{x^+ \in G(x,v)} B(x^+)\mid x].
\end{equation}}
\normalsize

As a result, the sufficient condition for reaching  $K_{\Delta}^C$ starting from $x$ with probability less than $\rho$ is given as
    \begin{equation} \label{eq:one_inv1}
        \Eeal[\sup_{x^+ \in G(x,v)} B(x^+)\mid x] \leq \rho \Delta.
    \end{equation}
For computing the bound of probability of reaching a the set $K_{\Delta}^C$ within $N$ time steps, we define sequences of functions, as follows.
Let $\lambda_{k}^B: \calX \times \calV^{k} \to \Real_{\geq 0}$ be defined recursively as follows:
\begin{equation}\label{eq:barrier-succ-funs}
\left.\begin{array}{rl}
\lambda_0^B(x) & = B(x) \\
\lambda_1^B(x, v_0) & = \displaystyle\sup_{x_1 \in G(x, v_0)} \lambda_0^B(x_1) \\ 
\vdots \\ 
\lambda_{k+1}^B(x, v_0, \ldots, v_k) & = \displaystyle\sup_{x_1 \in G(x, v_0)} \lambda_k^B(x_1, v_1, \ldots, v_k) \\ 
\end{array} \right\} 
\end{equation}
By Assumption \ref{ass:standing} that the set $G(x,v)$ for each $(x,v)$ is compact, then $\lambda_{k}$ is well-defined for each $k$.


\begin{lem}\label{lem:prob-k-steps}
Consider the system $\Sigma$ in \eqref{eq:sys} such that Assumptions \ref{ass:standing}--\ref{ass:iid} hold. Let $B:\calX\to \Real_{\geq 0}$ be a upper semicontinuous function with $\Delta$-sublevel set $K_\Delta$ is defined in \eqref{eq:k}.
For any $k \in \Neal$, probability of reaching $K_{\Delta}^C$ starting from $x \in \calX$ in \textbf{exactly} $k+1$th step is bounded from above by 
$\frac{1}{\Delta} \Eeal[\lambda_{k+1}^B(x, v_0, \ldots, v_k)\mid x]$.
\end{lem}
\begin{pf}
Proof with Markov's inequality.
\end{pf}

We denote $\ell_{k+1}^B(x) = \Eeal[\lambda_{k+1}^B(x, v_0, \ldots, v_k)\mid x]$.
\begin{lem} \label{lm:hb_suff}
    Consider the system $\Sigma$ in \eqref{eq:sys} such that Assumptions \ref{ass:standing}--\ref{ass:iid} hold. Let $B:\calX\to \Real_{\geq 0}$ be a upper semicontinuous function with $\Delta$-sublevel set $K_{\Delta}$ is defined in \eqref{eq:k}. The following is a sufficient condition for reaching the set $K_{\Delta}^C = \calX \setminus K_{\Delta}$ within $N$ steps with probability less than $\rho$ starting from $x \in K_{\Delta}$:
    \begin{equation}\label{eq:suff-cond-hb}
        \sum_{i=1}^N \ell_i(x) \leq \rho \Delta.
    \end{equation}
\end{lem}
\begin{pf}
Suppose~\eqref{eq:suff-cond-hb} holds.
    The probability of reaching the $K_{\Delta}^C:=\calX\backslash K_{\Delta}$ within $N$ steps is bounded by 
    Boole's inequality by the sum of the probabilities of reaching in precisely $j$ steps for $j$ 
    ranging from $0$ to $N$.
Trivially, since $x \in K_{\Delta}$, the probability of reaching $K_{\Delta}^C$ in $0$ steps is $0$. Using Lemma~\ref{lem:prob-k-steps}, 
we obtain the upper bound $\frac{1}{\Delta} \sum_{i=1}^N \ell_i^B(x) \leq \frac{1}{\Delta}  \Delta \rho \leq \rho$.

  \end{pf}
\begin{exmp} \label{ex:running_1}
Consider the following system
\begin{equation}
    x^+\in [x-v, x+v] 
\end{equation}
where $x\in \Real_{\geq 0}$, and $v\sim \calU[0,0.2]$. Let $B(x)=x$, then $\lambda_0^B(x) = x$, and 
$$\lambda_1^B(x, v_0) = \sup_{x_1 \in [x-v_0, x+v_0]} x_1 = x + v_0 \,.$$
Likewise, we can establish that 
\small{$$\lambda_{k}^B(x, v_0, \ldots, v_{k-1}) =  \displaystyle\sup_{x_1 \in [x-v_0, x+v_0]}  x_1 + \sum_{j=1}^{k-1} v_j = x + \sum_{j=0}^{k-1} v_j.$$}
Therefore, $\ell_{k}^B(x) = x + \displaystyle\sum_{j=0}^{k-1} \Eeal[v_j]= x + k\Eeal[v]$.
Hence, $\displaystyle\sum_{i=1}^N \ell_i^B(x) = Nx + \displaystyle\sum_{i=1}^N i \Eeal[v]= Nx + \frac{N(N+1)}{2 } \Eeal[v]$.
Considering that $\Eeal[v]=0.1$, we obtain that
$$\sum_{i=1}^N \ell_i^B(x) = Nx + \frac{N(N+1)}{20}. $$
Consequently, by Lemma \ref{lm:hb_suff}, choosing an initial state that satisfies
\begin{equation*}
    x \leq  \frac{ \rho \Delta}{N} - \frac{N+1}{20}  
\end{equation*}s
is sufficient to ensure that probability of reaching the set $\{x\in \Real_{\geq 0}\mid x>\Delta\}$ within \( N \) time steps is less than \( \rho \). 
Thus, for $N=4$, $\rho=0.25$, and $\Delta=20$, the condition is  $x \leq 1$.
\end{exmp}

\begin{rem}
One of the key difficulties of working with a barrier function lies in that of computing $\ell_{k+1}^B$ for $i \in \Neal$. 
In theory, this requires us to compute a function $\lambda_{k+1}^B$ that involves the random choices $v_0, \ldots, v_k$
and computing its expectation. In what follows, we present two key improvements: (a) In Section~\ref{sec:concave-barrier},
we will show that when the barrier function $B$ is concave and the set valued map $G$ is convex, we can work with 
$\Eeal(v_i)$ rather than $v_i$ itself; (b) In Section~\ref{sec:supermartingale}, we will impose a super-martingale 
condition on $B$ that will allow us to obtain bounds for arbitrary time horizon $k$. 
\end{rem}
\section{Concave Barrier Function}\label{sec:concave-barrier}
In this section, we consider barrier functions with extra concavity property. We show that we obtain better results for probabilistic safety.

\begin{defn}
    [Convex Set-valued Map]\label{def:convex_setvalued}
    A set-valued map $F$ is convex if its graph is convex. Therefore, $F$ is convex iff for each $x_1,x_2 \in \dom F$, and $\theta \in [0,1]$
    \begin{equation}
        \theta F(x_1) + (1-\theta) F(x_2) \subseteq F\big(\theta x_1 + (1-\theta)x_2\big).
    \end{equation}
    (See \cite[Definition $2.5.1$ and Lemma $2.5.2$ ]{Aubin:1991:VT:120830})
\end{defn}

In the following lemma, we aim to show that if $B$ is a concave function, and $F$ is a convex set-valued map, then under specific conditions, the marginal function inherits the concavity of $B$.
\begin{lem} \label{lm:concavity}
Consider a convex set-valued map $F:\calX \rightrightarrows \calX$ such that $F$ is upper semicontinuous and $F(x)$ for each $x$ is compact. 
Let $B : \calX \rightarrow \mathbb{R}_{\geq 0}$ be a nonnegative concave upper semicontinuous function. Suppose the function $f:\calX \to \mathbb{R}_{\geq 0}$ is defined as
\begin{equation} \label{eq:f_b_F}
    f(x):= \sup_{y \in F(x)} B(y).
\end{equation}
Then, the following conditions hold:
\begin{enumerate} [label={C\ref{lm:concavity}.\arabic*},leftmargin=*]
\item \label{item:ff1} $f$ is well-defined, nonnegative and upper semicontinuous.
\item \label{item:ff2} $f$ is concave.
\end{enumerate}
\end{lem}
\begin{pf}
To prove \ref{item:ff1}, since $F(x)$ is compact for each $x \in \calX$, the function $f$ is well-defined. As $f$ is defined by taking maximum of nonnegative function $B$ in a compact set, $f$ is itself nonnegative. Considering that $F$ is upper semicontinuous and has compact values, and the fact that $B$ is also upper semicontinuous, Proposition $2.9$ \cite{freeman2008robust} implies that $f$ is upper semicontinuous.

To prove \ref{item:ff2}, let $x_1,x_2 \in \calX$, and $\lambda \in [0,1]$, then with the assumption that $F$ is convex set-valued map, we have
\begin{equation}\label{eq:lmconcave}
\begin{aligned}
    f &(\theta x_1 + (1-\theta) x_2) = \sup_{y \in F(\theta x_1+ (1-\theta) x_2)} B(y)\\
    &\geq 
     \sup_{y \in \theta F(x_1)+ (1-\theta) F(x_2)} B(y)\\
     &= \sup_{z_1 \in F(x_1)} \sup_{z_2 \in F(x_2)} B(\theta z_1 + (1-\theta) z_2)\\
     &\geq
     \sup_{z_1 \in F(x_1)} \sup_{z_2 \in F(x_2)} \theta B(z_1) + (1-\theta) B(z_2)\\
     &= 
      \theta \sup_{z_1 \in F(x_1)}  B(z_1) + (1-\theta) \sup_{z_2 \in F(x_2)} B(z_2)\\
      &=\theta f(x_1) + (1-\theta) f(x_2).
\end{aligned}
\end{equation} 
\end{pf}

In our setting, we are interested on concavity of $\lambda_{i+1}^B$ defined in Eq.~\eqref{eq:barrier-succ-funs} with respect to random variable $v$. However, even though for concavity $\lambda_{1}^B$, it requires that $G$ be convex set-valued with respect to $v$, but for concavity of $\lambda_{i+1}^B$ for $i\geq 1$, the convexity of $G$ with respect to both $x$ and $v$ is required.
 \begin{assum} \label{ass:Gconvex}
 $G$ is convex set-valued map.
\end{assum}

The lemma below shows that if $B$ is concave, under Assumption \ref{ass:Gconvex}, the functions $\lambda^B$ are also concave.

\begin{lem} \label{lm:concave}
Consider system $\Sigma$ in \eqref{eq:sys} such that Assumptions \ref{ass:standing}, \ref{ass:iid}, and \ref{ass:Gconvex} hold.
Let $B : \calX \rightarrow \mathbb{R}_{\geq 0}$ be a concave upper semicontinuous function, and let $\lambda^B$ be functions defined in Eq.~\eqref{eq:barrier-succ-funs}. Then, $\lambda_{k}^B$ is concave for each $k\geq 0$.
\end{lem}
\begin{pf}
Using Assumption \ref{ass:Gconvex}, and by Lemma \ref{lm:concavity} we have that $\lambda_1^B$ is concave. For $i\geq 1$, we assume that $\lambda_i^B$ is concave, and prove that $\lambda_{i+1}^B$ is concave by Lemma \ref{lm:concavity}.
\end{pf}

The next lemma illustrates the condition for $N$-step reachability of $K_{\Delta}^C$ when the barrier function is concave.
\begin{lem} \label{lm:concave_res}
    Consider system $\Sigma$ in \eqref{eq:sys} such that Assumptions \ref{ass:standing}, \ref{ass:iid}, and \ref{ass:Gconvex} hold.
Let $B : \calX \rightarrow \mathbb{R}_{\geq 0}$ be a concave upper semicontinuous function, with $\Delta$-sublevel set $K_{\Delta}$. Let $\lambda^B$ be as defined in ~\eqref{eq:barrier-succ-funs}.
Then, the set $K_{\Delta}^C$ is reachable within $N$-step from $x \in K_{\Delta}$ with probability less than $\rho$, if 
\begin{equation} \label{eq:nstep_lambda}
\Eeal[\sup_{k \in [1, N]} \lambda_k^B(x,v_0,\cdots,v_{k-1})\mid x]\leq \rho \Delta.
\end{equation}
\end{lem}
\begin{rem}
When there exists $j \in [1,N]$ such that 
\begin{equation*}
    \sup_{k \in [1,N]} \lambda_k^B(x,v_0,\cdots,v_{k-1}) =  \lambda_j^B(x,v_0,\cdots,v_{j-1}),
\end{equation*}
using Lemma \ref{lm:concave_res}, we have that 
\begin{align*}
    \Eeal[\sup_{k \in [1, N]} \lambda_k^B(x,v_0,&\cdots,v_{k-1})\mid x]\\&= \lambda_j^B(x,\Eeal[v_0],\cdots,\Eeal[v_{j-1}]).
\end{align*}
Therefore, the condition \eqref{eq:nstep_lambda} simplifies to 
$$\lambda_j^B(x,\Eeal[v_0],\cdots,\Eeal[v_{j-1}]) \leq  \rho \Delta.$$
\end{rem}
\begin{exmp} \label{ex:running_2}
Consider Example \ref{ex:running_1}. The set-valued map dynamic is convex with respect to $v$, and the barrier function $B(x)=x$ is linear, thus concave. We have $\lambda_0=x$, and 
$\lambda_1(x,v_0) = \sup [x-v_0,x+v_0]= x+v_0$.
Also, 
$$\lambda_k^B(x, v_0, \cdots,v_{k-1})= x+\sum_{j=0}^{k-1}v_j.$$
Since $v_k$ is nonnegative, we have,
\begin{align*}
    \sup_{k\in [1,N]}  \lambda_k^B&(x, v_0,\cdots,v_{k-1})\\&=\lambda_N^B(x, v_0,\cdots,v_{N-1}) = x+\sum_{j=0}^{N-1}v_j 
\end{align*}
Using $\Eeal[v_i ]=0.1$ for each $i\geq 0$, we obtain that
$$\Eeal[\sup_{k\in [1,N]} \lambda_k^B(x, v_0,\cdots,v_{k-1})\mid x] =x+\frac{N}{10}.$$
Then, using Lemma \ref{lm:concave_res}, a sufficient initial condition for the probability of reaching the set $\{x\in \Real_{\geq 0}\mid x>\Delta\}$ within \( N \) time steps is less than \( \rho \) is given as $x\leq  \rho\Delta - \frac{N}{10}$.
Thus, for $N=4$, $\rho=0.25$, and $\Delta=20$, the condition is  $x\leq 4.6$.
\end{exmp}
\begin{rem}
    Let $w:=(\phi_1(v), \cdots, \phi_r(v))$, and suppose that the map $G(x,v)$ can be written as $\hat{G}(x, w)$, and  is a convex set-valued map over $w$, then the function $\lambda_k^B$s are still concave with respect to $w$ provided that $B$ is a concave function. This allows us to improve the applicability of Lemma~\ref{lm:concave_res} by considering ``higher moments'' of the random variables. 
\end{rem}
\begin{exmp} \label{ex:running_3}
Consider the following dynamics
\begin{equation}
    x^+ \in G(x,v)=[x+v-v^2, x+v+v^2].
\end{equation}
We can define $\zeta:=\begin{pmatrix}
    v \\ w
\end{pmatrix}$, where $w=v^2$. Then, replacing $v$ with $\zeta$, we have $G(x,\zeta)= [x+v-w,x+v+w]$, which is convex with respect to $\zeta$. Let $B(x)=x$. Then, we have
$$\lambda_k^B(x, \zeta_0, \cdots,\zeta_{k-1})= x+\sum_{j=0}^{k-1}(v_j+w_j).$$
Therefore,
\begin{align*}
    &\sup_{k\in [1,N]}  \lambda_k^B(x, \zeta_0,\cdots,\zeta_{k-1})\\&=\lambda_N^B(x, \zeta_0,\cdots,\zeta_{N-1}) = x+\sum_{j=0}^{N-1} (v_j +w_j). 
\end{align*}
Then, since $v\sim \calU[0,0.2]$, we have that $\Eeal[v_i ]=0.1$ for each $i\geq 0$, also $\Eeal[w_i ]=\Eeal[v_i^2 ]=\frac{1}{75}$.
Therefore,
$$\Eeal[\sup_{k\in [1,N]} \lambda_k^B(x, \zeta_0,\cdots,\zeta_{k-1})\mid x]=x+\frac{17N}{150}.$$
Then, using Lemma \ref{lm:concave_res}, a sufficient condition for bounding the probability of reaching the set $\{x\in \Real_{\geq 0}\mid x>\Delta\}$ within \( N \) time steps is less than \( \rho \) is given as $x\leq  \rho \Delta - \frac{17N}{150}$.
Thus, for $N=4$, $\rho=0.25$, and $
\Delta=20$, the condition is  $x \leq 4.547$.
\end{exmp}
In the next section, we analyze the safety properties when the barrier function satisfies the supermartingale condition.
\section{Concave Supermartingale Barrier Functions}\label{sec:supermartingale}
 A stochastic process is called a supermartingale if its expected value is nonincreasing with respect to the time. Here, we define the supermartingale barrier function for the system $\Sigma$ in \eqref{eq:sys}.

 \begin{defn}
[Barrier Function Candidate] A  nonnegative upper semicontinuous function $B:\calX \to \mathbb{R}_{\geq 0}$ is said to be a barrier function candidate with respect to initial and unsafe sets $(X_o, X_u)$, if  
\begin{equation}
\label{eq:candid}
\begin{aligned} 
 & B(x)  \leq \delta \qquad \forall x \in  X_o\\
 &  B(x)  \geq \Delta  \qquad \forall x \in X_u, \\
\end{aligned}
\end{equation}
where $\Delta > \delta\geq 0$.
\end{defn}
\begin{defn} 
    [Supermartingale Barrier Function] \label{def:supmB}
    Given the system $\Sigma$ in \eqref{eq:sys}, a barrier function candidate $B:\calX\to \Real_{\geq 0}$ is called supermartingale if for each $x \in \calX \backslash X_u$
 \begin{equation} \label{eq:supBar}
        \Eeal[\sup_{x^{+}\in G(x,v)} B(x^{+})\mid x] \leq B(x).
    \end{equation}
\end{defn}
To ensure safety for an unbounded time horizon, we seek a bound for the expectation of $\lambda_{k+1}^B$ functions. Since there is supremum function in defining $\lambda_{k+1}^B$ with respect to $\lambda_{k}$, we cannot extend bound of supermartingale to $\lambda_{k+1}^B$ for $k\geq 1$. However, when $\lambda_{k+1}^B$ is concave, the supermartingale condition can be extended.

The following results, we specify the condition of safety when $B$ is concave and the set-valued map is convex (see Definition \ref{def:convex_setvalued}).
\begin{lem} \label{lm:B_sup}
    Consider system $\Sigma$ in \eqref{eq:sys} such that Assumptions \ref{ass:standing}, \ref{ass:iid}, and \ref{ass:Gconvex} hold. Let $B:\calX \to \Real_{\geq 0}$ be a concave barrier function such that for each $x \in \calX \backslash X_u$, we have
    \begin{equation} \label{eq:supm_def}
        \sup_{x^{+}\in G(x,\Eeal[v])} B(x^{+})\leq B(x).
    \end{equation}
    Then, $B$ is a supermartingale barrier function.
\end{lem}
\begin{pf}
     The function $B$ is a concave function, and using Assumption \ref{ass:Gconvex}, $G$ is convex, therefore Lemma \ref{lm:concavity} implies that $\lambda_1^B(x,v): =\displaystyle\sup_{x^{+}\in G(x,v)} B(x^{+})$ is concave function. Therefore, we have
     $\Eeal[\lambda_1^B(x,v)\mid x] \leq \lambda_1^B(x,\Eeal[v])$.
     Using the condition \eqref{eq:supm_def}, we obtain that $\lambda_1^B(x,\Eeal[v]) \leq B(x)$. Consequently, we conclude that $\Eeal[\sup_{x^{+}\in G(x,v)} B(x^{+})\mid x] \leq B(x)$.    
\end{pf}

\begin{lem} \label{lm:h_B_sup}
    Consider system $\Sigma$ in \eqref{eq:sys} such that Assumptions \ref{ass:standing}, \ref{ass:iid}, and \ref{ass:Gconvex} hold. Let $B$ be a concave barrier function such that the condition \eqref{eq:supm_def} is satisfied.
   Then, for $k\geq 0$ and $x \in \calX\backslash X_u$, we have 
    \begin{equation}
        \label{lm:h_B_sup:item3}
        \lambda_k^B(x,\Eeal[v_0],\cdots,\Eeal[v_{k-1}]) \leq B(x).
    \end{equation}
\end{lem}
\begin{pf}
By induction.
For $k=0$ it is trivial. For $k=1$, it holds since it is the same condition as \eqref{eq:supm_def}.
Next, suppose the condition holds for $k$ and we show that it holds for $k+1$. 
Using concavity of $\lambda_i^B$, we have
\begin{align*}
\lambda_{k+1}^B(x,\Eeal[v_0],&\cdots,\Eeal[v_{k}])\\
    &= \sup_{x^{+} \in G(x,\Eeal[v_0])} \lambda_{k}^B(x^{+},\Eeal[v_1],\cdots,\Eeal[v_{k}]) \\
    &\leq \sup_{x^{+} \in G(x,\Eeal[v_0])} B(x^+)
    \\
    & \leq B(x).
\end{align*}
\end{pf}
\begin{defn} [Supermartingale]
   Given a probability space \((\Omega, \mathcal{F}, \mathbb{P})\), the stochastic process \(\{X_t\}_{t}\) is called a \emph{supermartingale} with respect to a filtration \(\calM = \{\calM_t\}_{t}\) iff the following conditions hold:
\begin{itemize}
    \item \(X_t\) is \(\calM_t\)-measurable for all \(t\),
    \item \(\Eeal[\vert X_t \vert] < \infty\) for all \(t\),
    \item \(\mathbb{E}[X_t \mid \mathcal{M}_s] \leq X_s\) for all $t,s$, \(t \geq s\)\,.
\end{itemize}
\end{defn}
Let $Z_k^x = \lambda_k^B(x,v_0,\cdots,v_{k-1})$. In the following lemma, we show that $\{Z_k\}$ has supermartingale property.
\begin{lem} \label{lm:superm_L}
     Consider system $\Sigma$ in \eqref{eq:sys} such that Assumptions \ref{ass:standing}, \ref{ass:iid}, and \ref{ass:Gconvex} hold. Let $B$ be a concave barrier function such that the condition \eqref{eq:supm_def} is satisfied.
   Then, for $k\geq 0$ and $x \in \calX\backslash X_u$, the sequence $\{Z_k^x\}$, where $Z_k^x = \lambda_k^B(x,v_0,\cdots,v_{k-1})$ is supermartingale.
\end{lem}
\begin{pf}
    Let $k \geq j$. We have that
    \begin{align*}
        &\Eeal[Z_k \mid Z_j] = \\
        &\Eeal[ \lambda_k^B(x,v_0,\cdots,v_{k-1})| v_0=w_0,\cdots, v_{j-1}=w_{j-1})]\\
        & \leq \lambda_k^B(x,w_0,\cdots,w_{j-1},\Eeal[v_{j}], \cdots, \Eeal[v_{k-1}])\\
        & = \sup_{y_1\in G(x,w_0)} \cdots\sup_{y_j\in G(y_{j-1},w_{j-1})} \lambda_{k-j}^B(y_j, \Eeal[v_{j}], \cdots, \Eeal[v_{k-1}])\\
        &\leq \sup_{y_1\in G(x,w_0)} \cdots\sup_{y_j\in G(y_{j-1},w_{j-1})} B(y_j)\\
        &= \lambda_j^B(x, w_0, \cdots, w_{j-1})= Z_j.
    \end{align*}
\end{pf}
\begin{rem}
    Using \cite[Definition $10.1.2$]{oksendal2013stochastic} and Lemma \ref{lm:h_B_sup}, we can conclude that $\{\lambda_k^B\}$ are called \emph{supermeanvalued} with respect to $v$. Also, \cite[page $216$]{oksendal2013stochastic} indicates that the sequence is a supermartingale.
\end{rem}

\begin{thm}
    Consider system $\Sigma$ in \eqref{eq:sys} such that Assumptions \ref{ass:standing}, \ref{ass:iid}, and \ref{ass:Gconvex} hold. Let $B$ be a concave nonnegative barrier function candidate with respect to $(X_0,X_u)$ such that the condition \eqref{eq:supm_def} is satisfied. Then, the probability of reaching $X_u$ starting from $X_0$ is less than $\frac{\delta}{\Delta}$.
\end{thm}
\begin{pf}
    Let $N \in \overline{\Neal}$, and $x_0 \in X_0$. The probability that the solution to the system $\Sigma$ hits the unsafe set $X_u$ in less than $N$ time steps is given by
    \begin{align*}
        \Peal(\exists k & \in [1,N]\,, \lambda_k^B(x_0,v_0,\cdots,v_{k-1})\geq \Delta \mid x_0) \\
        &= \Peal(\sup_{k\in [1,N]} \lambda_k^B(x_0,v_0,\cdots,v_{k-1}) \geq \Delta \mid x_0)
    \end{align*}
    Next, Lemma \ref{lm:superm_L} implies that $\{\lambda_k^B\}$ is a supermartingale. Therefore, using supermartingale inequality, and Lemma \ref{lm:concave} we have 
    \begin{align*}
        \Peal(\sup_{k\in [1,N]} \lambda_k^B(&x_0,v_0,\cdots,v_{k-1}) \geq \Delta \mid x_0) \\
        & 
        \leq \frac{1}{\Delta}\lambda_0(x_0)= \frac{1}{\Delta} B(x_0).
    \end{align*}
Therefore, the probability of reaching $X_u$ from $X_0$ is less than $\displaystyle\sup_{x_0\in X_0} \frac{1}{\Delta}B(x_0)= \frac{\delta}{\Delta}$.
\end{pf}

\section{Numerical Examples}
We illustrate the results with some examples.
\begin{exmp} \label{ex:lin}
    Consider the following system
    \begin{equation}
        x^+ \in \left\{ \left(\gamma A_1 + (1-\gamma) A_2\right)x + b v\ |\ \gamma \in [0, 1] \right\},
    \end{equation}
    where $A_1,A_2 \in \Real^{n\times n}$ and $b \in \Real^n$, with $v$ denoting the random variable. 
    The sate space is denoted by $\calX\subseteq \Real^n $, and the unsafe set is given by $X_u \subseteq \calX$.
    We consider a linear barrier function $B(x)= c_0+c^\top x$, defined over $\calX$. Since $B$ is linear, it is also concave. The system dynamic is linear with respect to $(x,v)$, so it is convex.
    Let $\Tilde{A}(\gamma) := \gamma A_1 + (1-\gamma) A_2$, then we have
    \begin{align*}
        \lambda_0(x) &= c_0+c^\top x \\
        \lambda_1(x,v_0) &= \sup_{\gamma_0 \in [0,1]} c_0+c^\top bv_0 + c^\top \Tilde{A}(\gamma_0)x
    \end{align*}
    \begin{align*}
        &\lambda_i(x,v_0 ,\cdots,v_{i-1}) = \sup_{\gamma_{0} \in [0,1]}\cdots \sup_{\gamma_{i-1} \in [0,1]} c_0+c^\top bv_{i-1}\\
        & + c^\top \Tilde{A}(\gamma_{i-1})bv_{i-2}+\cdots+ c^\top \Tilde{A}(\gamma_{i-1})\cdots \Tilde{A}(\gamma_{1})bv_{0}\\
        &+
        c^\top \Tilde{A}(\gamma_{i-1})\cdots \Tilde{A}(\gamma_{0})x\,.
    \end{align*}
     
    Suppose $X_u$ is subset of $\Delta$-superlevel set of $B$:
    \[ X_u = \{ x \in \calX \ |\ B(x) \geq \Delta \} \,,\]
    then the following condition specifies the safety of $x$ for $N$ time step with probability $1-\rho$:
    \begin{align}
        \Eeal[\sup_{i\in [1,N] } \lambda_i(x,v_0 ,\cdots,v_{i-1})\mid x] \leq \rho \Delta \,.
    \end{align}
    Furthermore, the concave supermartingale condition is given as
     \begin{equation} \label{eq:lin_ex_sup}
     \lambda_1(x,\Eeal[v_0]) = \sup_{\gamma_0 \in [0,1]} c_0+c^\top b\Eeal[v_0] + c^\top \Tilde{A}(\gamma_0)x\leq c_0+c^\top x.
     \end{equation}
     Since the objective function in the optimization \eqref{eq:lin_ex_sup} is linear with respect to $\gamma$, so the optimal value is in the boundary of set $[0,1]$. Then, the supermartingale condition reduces to the following for $x \in \calX \backslash X_u$
     \begin{equation} \label{eq:lin_ex_sup2}
     \begin{aligned} 
         &c^\top \left((A_1-I)x+b \Eeal[v] \right)\leq 0\\
         &c^\top \left( (A_2-I)x+b \Eeal[v]\right)\leq 0.
     \end{aligned}
     \end{equation}
     Hence, if the supermartingale condition in \eqref{eq:lin_ex_sup2} holds, every initial set equals to $\delta$-sublevel set of $B$, where $\delta \leq \Delta$, is safe with probability of $1-\frac{\delta}{\Delta}$.

     Suppose $A_1=\begin{pmatrix}
         -0.2 & 1\\ 0 & 0.3
     \end{pmatrix}$ $A_1=\begin{pmatrix}
         0.1 & -1\\ 0 & 1
     \end{pmatrix}$, $b=\begin{pmatrix}
         0.1 \\ -0.1
     \end{pmatrix}$, and let 
     $\calX = [0,5]^2$, $X_u =\{x \in \calX\mid x_2 \geq 3\}$.
     Suppose $v\sim \calU[1,2]$. 
     Then, considering $B(x)= x_2$, we have that
     $$\lambda_1(x,\Eeal[v_0]) = \sup_{\gamma_0 \in [0,1]} 
     \begin{pmatrix}
         0 & 1
     \end{pmatrix} \left(1.5 b +\Tilde{A}(\gamma_0)x\right).$$
     Since the above optimization is linear in $\gamma_0$, therefore, the answer is either when $\gamma_0=0$ or $\gamma_0=1$. Thus,
     \begin{align*}
         \lambda_1&(x,\Eeal[v_0]) \\
         =&\sup \left(-0.15 + \begin{pmatrix}
         0&1
     \end{pmatrix} A_1 x,  -0.15+ \begin{pmatrix}
         0&1
     \end{pmatrix} A_2 x\right)\\
     =& \sup (-0.15+0.3x_2, -0.15+x_2) 
     \end{align*}
Then, we have $\lambda_1(x,\Eeal[v_0])=  -0.15+x_2$ for $x \in \calX$.
Since $\lambda_1(x,\Eeal[v_0])=  -0.15+x_2\leq B(x)=x_2$ for $x\in \calX$, therefore, $B$ is a supermartingale for the system. Since $B(x)\geq 3$ for $x \in X_u$, then the probability of safety each for $\delta$-sublevel set of $B$, i.e,  $X_0=\{x\in \calX \mid x_2\leq \delta\}$,  is given as $1-\frac{\delta}{3}$. Figure \ref{fig:linear} shows $200$ trajectories for $N=30$ from randomly generated initial points in $X_0 = [0,5]\times [0,1]$ to the system \eqref{ex:lin}.

 \begin{figure}[t]
\centering
  \includegraphics[scale=0.5]{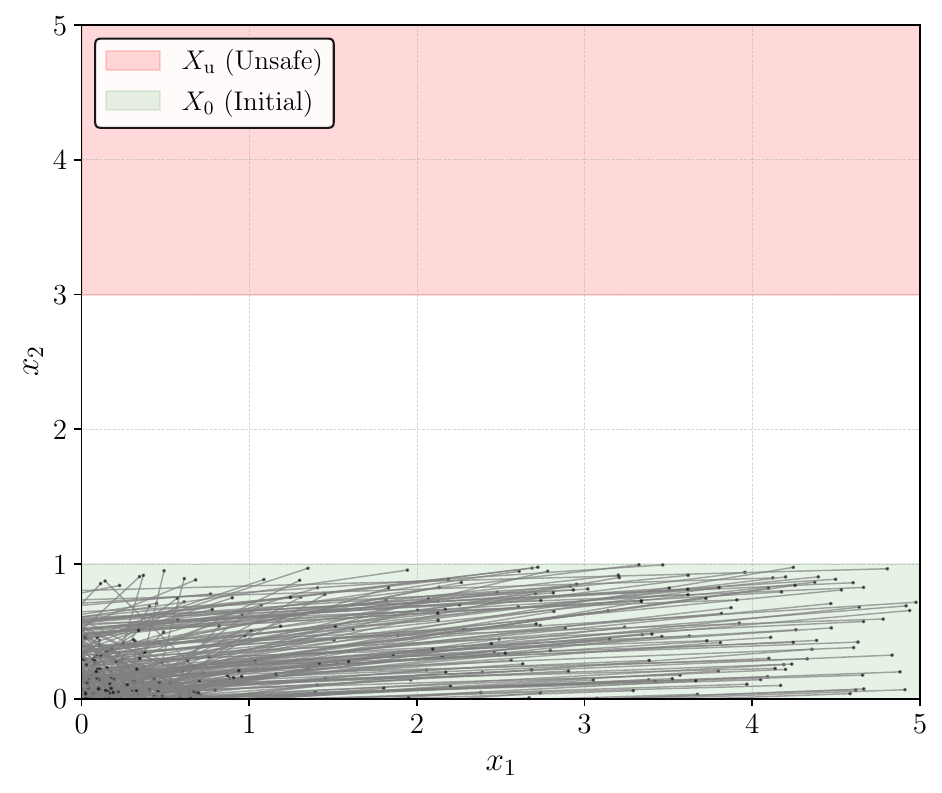}
  \caption{$200$ trajectories from randomly generated initial points in $X_0 = [0,5]\times [0,1]$ to the system \eqref{ex:lin}.}
  \label{fig:linear}
\end{figure}
\end{exmp}
\begin{exmp} \label{ex:matrix}
    Consider the following dynamics
    \begin{equation} \label{eq:sys_ex_mat}
        X^+ \in \left\{ U X U^\top+ M(\gamma) v\ | \ \gamma \in [a, b] \right\}, 
    \end{equation}
    where $X \in S^{n}_{++}$ is the state, $U$ is a unitary matrix, $M(\gamma) \in S^n$ is a linear matrix-valued map with respect to $\gamma$.  $v$ is a random variable.

    Suppose $\calX=S^{n}_{++}$, and the unsafe is defined as $X_u= \{X \in \calX\mid \Tr(X)\geq \Delta\}$. Let $B(x)= \Tr(X)$, then $B$ is positive on $\calX$.  Since the dynamics are linear with respect to $X$ and $v$, so it is convex. Also, the barrier function $B$ is linear so it is concave.
Then, we have
\begin{align*}
    \lambda_1(X,v_0) &= \sup_{\gamma_0 \in [a,b]} \Tr\big(U X U^\top+ M(\gamma_0) v_0 \big) \\
    &= \sup_{\gamma_0 \in [a,b]}  \Tr(X)+v_0\Tr(M(\gamma_0))\,.
\end{align*}
Also, $\lambda_i$ is given as
\begin{align*}
    \lambda_i(X,v_0,&\cdots,v_{i-1}) =\\
    &
    \sup_{\gamma_{0}  \in [a,b]}\cdots \sup_{\gamma_{i-1}  \in [a,b]} \Tr (X )+ \sum_{j=0}^{i-1} v_j \Tr(M(\gamma_j))\,. 
\end{align*}    

The condition specifies safety of $X$ for $N$ time step with probability $1-\rho$ is given as
    \begin{align}
        \Eeal[\sup_{i\in [1,N] } \lambda_i(X,v_0 ,\cdots,v_{i-1}) \mid x] \leq \rho \Delta .
    \end{align}
    
Suppose that $M(\gamma) = \gamma A_1 + (1 - \gamma) A_2$, where 
$A_1 = \begin{bmatrix}1 & 0.5 \\ 0.5 & 1\end{bmatrix}$ and $A_2 = \begin{bmatrix}0 & -1 \\ -1 & 0\end{bmatrix}$, and $\gamma\in [0,1] $.
The unitary matrix $U$  is expressed as a rotation by $\theta = \frac{\pi}{4}$, and $v \sim \calU[-0.2,0.2]$. Let $\Delta=1.2$, and $\delta=0.5$.

Therefore, $\Tr(M(\gamma_j))= 2\gamma_j$, and we have
$$\lambda_i(X,v_0,\cdots,v_{i-1}) = \Tr (X )+ 2\sum_{j=0}^{i-1} v_j \calI_{[0,0.2]}(v_j)\,.$$
As a result 
\begin{align*}
    \sup_{i\in [1,N]} \lambda_i(X,v_0,\cdots,v_{i-1})&= \lambda_N^B(X,v_0,\cdots,v_{N-1})\\
\end{align*}
Then, the condition that specifies safety of $X$ for $N$ time step with probability $1-\rho$ is given 
\begin{align*}
    \lambda_N(X,\Eeal[v_0],\cdots,\Eeal[v_{N-1}])\leq \rho \Delta.
\end{align*}
Therefore,
\begin{align*}
    \lambda_N(X,\Eeal[v_0],&\cdots,\Eeal[v_{N-1}]) =\\&\Tr (X )+ 2 \sum_{j=0}^{N-1} \Eeal[v_j] \calI_{[0,0.2]}(\Eeal[v_j])\\
    &= \Tr (X )
\end{align*}
Hence, the condition is $\Tr (X ) \leq \rho \Delta$.

    Furthermore, the supermartingale condition, namely,
     \begin{equation}
     \lambda_1(X,\Eeal[v_0]) =  \sup_{\gamma_0 \in [a,b]}  \Tr(X)+\Eeal[v_0]\Tr(M(\gamma_0)) \leq \Tr(X),
     \end{equation}
     holds on $\calX$ if 
     $$\displaystyle\sup_{\gamma_0 \in [a,b]} \Eeal[v_0]\Tr(M(\gamma_0)) \leq 0\,.$$
     Then, every initial set that is given as a $\delta$-sublevel set of $B$, where $\delta \leq \Delta$, is safe with probability of $1-\frac{\delta}{\Delta}$.

Figure \ref{fig:matrix} shows trace trajectories of $200$ randomly generated initial $X$ from $\Tr(X_0) < 0.3$ for $N=30$.
Since $E[v]=0$, then the barrier function $B$ is concave supermartingale, and the probability of reaching unsafe set, $\mathrm{trace}(X_t) > 1.2$, is given by $1-\frac{0.3}{1.2}= 0.25$. Figure \ref{fig:matrix} simulates  $200$ trajectories. We note that $12$ of these trajectories have entered the unsafe region.

     \begin{figure}[t]
\centering
  \includegraphics[scale=0.4]{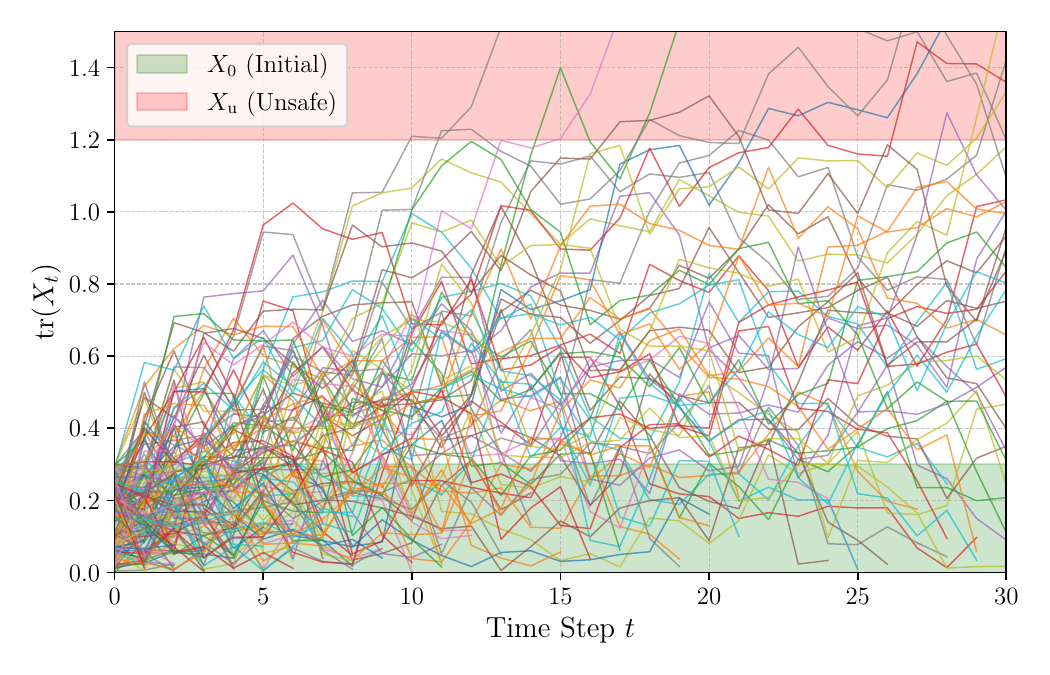}
  \caption{The figure shows the $200$ trajectories of trace from randomly generated initial points in $X_0 = \{X\in \calX\mid \Tr(X) \in (0,0.3]\}$ to the system \eqref{eq:sys_ex_mat} in Example \ref{ex:matrix}.}
  \label{fig:matrix}
\end{figure}
\end{exmp}

\section{Conclusion} \label{sec:conclusion}
In this work, we study the safety of stochastic difference inclusion systems through the use of barrier functions. The key contribution lies in a framework that allows for the easier verification of barrier functions by the use of concave barrier functions and convex set-valued maps for bounding probabilities of reaching unsafe sets over unbounded time horizons. In the future, we wish to study the application of these notions for other properties including stability and set-valued maps arising from problems with unknown distribution families defined by notions such as the Wasserstein metric. 

\appendix
\section{Measurability}
Here, we present some useful theorems on the measurability of compositions of set-valued maps from \cite{rockafellar2009variational} and provide the proofs of Lemmas \ref{lm:meas} and \ref{lm:meas_well}.

\subsection{Measurability of Set-valued Map}
\begin{thm} [Composition of Mappings] \label{thm:appendix_1413}
    \cite[Theorem~14.13(b)]{rockafellar2009variational}
    Let $S : T \rightrightarrows \mathbb{R}^n$ be a closed-valued and measurable set-valued map, and for each $t \in T$  
consider an outer semicontinuous mapping $M(t, \cdot) : \mathbb{R}^n \rightrightarrows \mathbb{R}^m$. Suppose the graphical mapping 
$t \mapsto \operatorname{gph} M(t, \cdot) \subset \mathbb{R}^n \times \mathbb{R}^m$
is measurable. Then the mapping  
$t \mapsto M(t, S(t))$ is measurable.
\end{thm}
\begin{cor}[Composition with Functions]
    \cite[Corollary $14.14$]{rockafellar2009variational}
    For each $t \in T$, consider an outer semicontinuous mapping  
$M(t, \cdot) : \mathbb{R}^n \rightrightarrows \mathbb{R}^m$,
and suppose the graphical mapping  
$
t \mapsto \operatorname{gph} M(t, \cdot) \subset \mathbb{R}^n \times \mathbb{R}^m
$
is measurable.  
Then, whenever $t \mapsto x(t)$ is measurable, the mapping  
$
t \mapsto M(t, x(t))
$ is closed-valued and measurable.
\end{cor}

\begin{pf}[Proof of Lemma \ref{lm:meas}] \label{pf:lm:meas}
 Under Assumption \ref{ass:standing}, the map $G$ is upper semicontinuous with respect to $(x,u)$ and has closed values. This implies that $G$ is also outer semicontinuous with respect to $x$.
 For $k=0$, under Assumption \ref{ass:standing}, \cite[Corollary $14.14$]{rockafellar2009variational} implies that $G(x,v_0)$ is measurable and closed.
 Suppose this is valued for $k$, we show it for $k+1$.
 Let $H_k := G(\cdots (G(G(x,v_0),v_1),v_2),\cdots),v_k)$, then $H_{k+1}=G(H_k,v_{k+1})$.
 Since $H_k$ is measurable ans closed, Theorem \ref{thm:appendix_1413} implies that $H_{k+1}$ is measurable. Also, the graph of $G$ is closed for every closed domain (see \cite[Proposition~1.4.8]{aubin2009set}), hence $H_{k+1}$ is closed.
\end{pf}
\subsection{Normal Integrand and Measurability Marginal Function}
\begin{defn}
    [normal integrands] \cite[Definition $14.27$]{rockafellar2009variational}
    A function $f:T \times \Real^n \to \bar{\Real}$ is called a \emph{normal integrand} if its epigraphical mapping \( S_f : T \to \mathbb{R}^n \times \mathbb{R} \) is closed-valued and measurable, where
    $S_f(t):=\{(x,\alpha) \Real^n \times \Real \mid f(t,x)\leq \alpha\}.$
\end{defn}

\begin{pf}[Proof of Lemma \ref{lm:meas_well}] \label{pf:lm:meas_well}
    Using  $\delta_C(x):=\begin{dcases}
    0 & x \in C\\
    -\infty & \text{otherwise}
\end{dcases}$, the function $\lambda$ can be written as $\lambda(x,v):=\displaystyle\sup_{y \in G(x,v)} h(y)=\displaystyle\sup_y \{h(y)+\delta_{G(x,v)}(y)\}$.
Let $\Bar{\lambda}(x,v):=-\lambda(x,v)=\inf_y \{-h(y)-\delta_{G(x,v)}(y)\}$.
Since $h$ is upper semicontinuous, then $-h$ is lower semicontinuous, then \cite[Example $14.31$]{rockafellar2009variational} implies that $-h$ is a \emph{normal integrand}. Next, since Lemma \ref{lm:meas} implies that $G(x,v)$ is measurable, then \cite[Example $14.32$]{rockafellar2009variational} implies that $-\delta_{G(x,v)}(y)= \begin{dcases}
    0 & x \in G(x,v)\\
    \infty & \text{otherwise}
\end{dcases}$ is normal integrand. Since the summation of normal integrands is normal integrand, \cite[Proposition $14.44(c)$]{rockafellar2009variational}, $f(x,v,y)=-h(y)-\delta_{G(x,v)}(y)$ is normal integrand. 
Since $h$ is upper semicontinuous, and $G$ is also upper semicontinuous with compact values, \cite[Proposition $2.9$]{freeman2008robust} implies that $\lambda$ is upper semicontinuous with respect to $x$. Then, $\Bar{\lambda}(x,v)$ is lower semicontinuous with respect to $x$. As a result, 
\cite[Proposition $14.47$]{rockafellar2009variational} implies that $\Bar{\lambda}$ is normal integrand. Therefore, it is measurable \cite[Proposition $14.34$]{rockafellar2009variational}.
Consequently, $\lambda$ is measurable.
Moreover, since $h$ is nonnegative, $\lambda$ is nonnegative. Therefore $I_\lambda[x]$ is well-defined in $[0,+\infty]$ for each $x \in \calX$.

 \end{pf}
\subsection{Measurability of $\lambda_{k+1}$}
Consider the set valued map $G:\calX \rightrightarrows \Real^n$. For finite $N \in \Neal$, let $F_N:\calX \rightrightarrows \Real^n$ be a set-valued map which is defined as 
\begin{equation} \label{eq:F_compose}
    F_N(x) := \underbrace{G \circ G \circ \cdots \circ G}_{N \text{ times}}(x).
\end{equation}
 
\begin{lem} \label{lm:compos_usc}
    Let $G:\calX \rightrightarrows \Real^n$ be a upper semicontinuous map which is locally bounded with closed values. For finite $N \in \Neal$, consider the set-valued map $F_N:\calX \rightrightarrows \Real^n$ in \eqref{eq:F_compose}. The map $F$ is upper semicontinuous map with compact values.
\end{lem}
\begin{pf}
    Using \cite[Proposition $5.52$ (b)]{rockafellar2009variational}, we conclude that $F_N$ is upper semicontinuous. 
    Since $G$ has compact values, then using \cite[Proposition $2.6$]{freeman2008robust}, for every compact set $K \subset \calX$, $G(K)$ is compact. Therefore, $F_N(x)$ is compact.
\end{pf}
\begin{pf}[Proof of Lemma \ref{lm:lambda_meas_k}] \label{pf:lm:lambda_meas_k}
Using that indicator function is upper semicontinuous for a closed set $S$, Lemma \ref{lm:meas_well} implies that $\lambda_1$ is measurable.
We can formulate the function $\lambda_{k+1}$ as
\begin{align} \label{lambda_GG}
    \lambda_{k+1}^S(x, v_0, \ldots, v_{k})= 
    \sup_{y_{k+1} \in G(\cdots (G(x,v_0),v1),\cdots,v_k)} \calI_S(y_{k+1}). 
\end{align}
Let $H(x,v_0,\cdots,v_k):=G(\cdots (G(x,v_0),v1),\cdots,v_k)$.
Using Lemma \ref{lm:compos_usc}, we conclude that $H$ is upper semicontinuous with respect to $x$. 
Also, Lemma \ref{lm:meas} implies that $H$ is measurable with closed values. Therefore, Lemma \ref{lm:meas_well} implies that $\lambda_{k+1}$ for $k \in \mathbb{N}$ is measurable.
\end{pf}

\bibliographystyle{plain}        
\bibliography{references}           



\end{document}